\documentclass[epic,11pt]{article}
\usepackage{amsfonts,amsmath,amscd,amssymb,amsthm}
\usepackage{stackrel}
\usepackage{fullpage}
\usepackage{xcolor}
\usepackage[all]{xy}
\usepackage[alphabetic,initials]{amsrefs}
\input xy%
\usepackage{fullpage}
\usepackage{amsrefs}
\makeindex

\def\convf{\hbox{\space \raise-2mm\hbox{$\textstyle      \bigotimes \atop \scriptstyle \omega$} \space}}
\def\0{{\bar 0}}
\def\1{{\bar 1}}

\def\Z{{\mathbb Z}}

\def\K{{\mathbb K}}
\def\E{{\mathbb E}}
\def\L{{\mathbb L}}
\newcommand{\ra}{\rightarrow}
\def\R{{\mathbb R}}

\def\I{{\mathbb I}}
\def\du{\stackrel{\cdot}{\bigcup}}
\newcommand{\br}{\begin{rem}}
\newcommand{\er}{\end{rem}}
\newcommand{\lp}{\mathtt{lp}}
\newcommand{\ttT}{\mathtt{T}}
\newcommand{\noi}{\noindent}
\newcommand{\itema}{\item[{{\rm(a)}}]}
\newcommand{\itemb}{\item[{{\rm(b)}}]}
\newcommand{\cF}{\mathcal{F}}
\def\LM{\operatorname{LM}}
\def\PM{\operatorname{PM}}

\def\LN{\operatorname{PM^*}}
\newcommand{\fg}{\mathfrak{g}}\newcommand{\fgl}{\mathfrak{gl}}

\newcommand{\ci}{\cdot} \newcommand{\ti}{\times}
\newcommand{\ey}{\end{eqnarray}}
\newcommand{\by}{\begin{eqnarray}}
\newcommand{\nn}{\nonumber}

\newcommand{\ec}{\end{corollary}}

\newcommand{\bpf}{\begin{proof}}
\newcommand{\epf}{\end{proof}}
\newcommand{\bt}{\begin{theorem}}
\newcommand{\et}{\end{theorem}}
\newcommand{\bi}{\begin{itemize}}
\newcommand{\ei}{\end{itemize}}
\newcommand{\bl}{\begin{lemma}}
\newcommand{\bp}{\begin{proposition}}
\newcommand{\be}{\begin{equation}}
\newcommand{\bc}{\begin{corollary}}
\newcommand{\bexs}{\begin{examples}}
\newcommand{\eexs}{\end{examples}}
\newcommand{\bexa}{\begin{example}}
\newcommand{\eexa}{\end{example}}

\newcommand{\el}{\end{lemma}}
\newcommand{\ep}{\end{proposition}}
\newcommand{\ee}{\end{equation}}
\newcommand{\cB}{\mathcal{B}}
\newcommand{\cL}{\mathcal{L}}
\newcommand{\cK}{\mathcal{K}}

\begin{document}
\title{Table of Contents}
\newtheorem{definition}{Definition}
\newtheorem{thm}{Theorem}[section]
\newtheorem{hyp}[thm]{Hypothesis}
 \newtheorem{hyps}[thm]{Hypotheses}
  \newtheorem{rems}[thm]{Remarks}

\newtheorem{conjecture}[thm]{Conjecture}
\newtheorem{theorem}[thm]{Theorem}
\newtheorem{theorem a}[thm]{Theorem A}
\newtheorem{example}[thm]{Example}
\newtheorem{examples}[thm]{Examples}
\newtheorem{corollary}[thm]{Corollary}
\newtheorem{rem}[thm]{Remark}
\newtheorem{lemma}[thm]{Lemma}
\newtheorem{sublemma}[thm]{Sublemma}
\newtheorem{cor}[thm]{Corollary}
\newtheorem{proposition}[thm]{Proposition}
\newtheorem{exs}[thm]{Examples}
\newtheorem{ex}[thm]{Example}
\newtheorem{exercise}[thm]{Exercise}
\numberwithin{equation}{section}%
\setcounter{part}{0}

\numberwithin{definition}{section}
\newcommand{\drar}{\rightarrow}
\newcommand{\lra}{\longrightarrow}
\newcommand{\llra}{\longleftrightarrow}


\newtheorem*{thm*}{Theorem}
\newtheorem{lem}[thm]{Lemma}
\newtheorem*{lem*}{Lemma}
\newtheorem*{prop*}{Proposition}
\newtheorem*{cor*}{Corollary}
\newtheorem{dfn}[thm]{Definition}
\newcommand{\bd}{\begin{definition}}
\newcommand{\ed}{\end{definition}}

\newtheorem*{defn*}{Definition}
\newtheorem{notadefn}[thm]{Notation and Definition}
\newtheorem*{notadefn*}{Notation and Definition}
\newtheorem{nota}[thm]{Notation}
\newtheorem*{nota*}{Notation}
\newtheorem{note}[thm]{Remark}
\newtheorem*{note*}{Remark}
\newtheorem*{notes*}{Remarks}
\newtheorem{hypo}[thm]{Hypothesis}
\newtheorem*{ex*}{Example}
\newtheorem{prob}[thm]{Problems}
\newtheorem{conj}[thm]{Conjecture}
\newcommand{\Bc}{\begin{center}}
\newcommand{\Ec}{\end{center}}
%
\title{ Catalan numbers and  a conjecture on the maximum composition length of a  Kac module}
\author{Ian M. Musson \\
University of Wisconsin-Milwaukee\\ email: {\tt
musson@uwm.edu}}

\maketitle
\begin{abstract} Let $f:\mathbb{Z}\lra \{ \ti, \ci\}$ be a function such that $f(a) = \ci$ for all except finitely for 
 many  $a \in \mathbb{Z}$.   We define a set $\flat f$ of non-intersecting arc (or cap) diagrams satisfying certain conditions determined by $f$. Then we give a  recursive method for  enumeration of $\flat f$ which  recalls the Fundamental Recurrence for Catalan numbers. The motivation comes from the problem of enumeration of the 
 composition factors of a  Kac module with maximum degree of atypicality for the Lie superalgebra $\fg=\fgl(r|r)$. In particular we show the  maximum number of 
 composition factors is a Catalan number, as conjectured in \cite{MSe}. 
\end{abstract}
\Bc Keywords:  Catalan numbers, Kac modules \Ec
\section{Introduction} 
Many problems in representation theory involve finding the composition factors of some natural module $M$, such as a (parabolic) Verma module, Specht module or Kac module.  If $M$ is multiplicity-free and the composition factors are indexed by some combinatorially defined set $S$, the problem might be considered solved.  However we then have the  combinatorial problem of how to enumerate $S$.
\\ \\
For Kac modules, we give a {\it Catalan-type recursion} for the set of composition factors.  To do this we use legal moves as defined in \cite{MSe} and  tally functions which, though completely elementary, give a new way of looking at these moves. 

\subsection{Kac modules} 
 Let $X^+ =X^+(m|n)$ denote the set of dominant
integral weights for 
the complex Lie superalgebra $\fg = \fgl(m|n)$. 
In  \cite{Kac2} Kac introduced a certain finite dimensional
highest weight module $K(f)$ 
with highest weight   $f \in X^+$, whose character is given by an analog of  
the Weyl character formula.  The Kac module 
$K(f)$ has a unique simple factor module $L(f)$ and  any composition factor of $K(f)$ is isomorphic to some  $L(g)$ with $g \in X^+$.  This raises the problem of  finding the  
multiplicity of $L(g)$  as a composition factor of $K(f).$ 
Combinatorial solutions to this problem 
 were given  in  \cite{S2}  and
\cite{Br}, using completely different methods. 
\\ \\
\noi  Let $\cF$ be the category of finite
dimensional $\fg$-modules which are weight modules for a Cartan subalgebra.  
In Section \ref{ab} we recall the notion of matching cap  and weight diagrams.  
  If $f \in X^+,$ we set 
\be \label{VGL} \flat (f) = \{g \in
X^+|D_{wt}(f)\; \mbox{matches the cap diagram} \;D_{cap}(g)\}.\ee
The main result of \cite{Br} can be expressed as follows, see \cite{MSe}.

\bt \label{brnt} Each Kac module is multiplicity free and in the Grothendieck group $\K[\cF]$ we have
\be \label{KGL} K(f) = \sum_{g\in \flat (f)} L(g).\ee \et \noi 
Thus $\flat (f)$ is the set we wish to enumerate.  We note another occurence of this set in representation  theory.    By \cite{BS4} Theorem 1.1 if $\cB$ is any block of $\cF$, there is an equivalence of categories $\E$ from $\cB$ to a block 
of a Khovanov diagram algebra such that $\E  K(f) =\cK(f)$ is a standard module and $\E L(f) =\cL(f)$ is simple. 
As an immediate consequence, we have in a suitable  
Grothendieck group, (see also  
\cite{BS4} Theorem 2.1.)
  \be \label{SGL} \cK(f) = \sum_{g\in \flat (f)} \cL(g).\ee 
By an equivalence of categories of Serganova \cite{S3}, \cite{GS} any block of the category $\cF$ is equivalent to the most atypical block  $\cB$ for 
$\fg=\fgl(r|r)$ and from now on we  work only with this block. This reduction simplifies the notation that we introduce in Subsection \ref{ab}.  Without this reduction the maximum composition length of a Kac module in  any block will be the Catalan number $C_{r+1},$ 
where $r$  is the 
 degree of atypicality of the block. Denote the  set of highest weights in $\cB$ by $X(r)$.

\subsection{Catalan numbers}  \begin{definition} {\rm Let $C_0=1$.  For $n\ge 1$ the $n^{th}$ {\it Catalan number} $C_n$ can be defined as the number of   ways of connecting $2n$ points lying on a horizontal line by $n$ nonintersecting arcs (or caps), each arc connecting two of the points and lying above the points, \cite{St16} Bijective Exercise 61.}\end{definition} \noi 
 Hopefully  the diagrams in  this paper will make any unexplained terminology clear.  Here are the $C_3=5$ ways of connecting 6 points. 
\\
\setlength\unitlength{0.48mm}
\begin{picture}(30,20)
\put(0,2){\line(1,0){50}}
\multiput(0,2)(10,0){6}{\circle*{1}}
\put(5,2){\oval(10,8)[t]}
\put(25,2){\oval(10,8)[t]}
\put(45,2){\oval(10,8)[t]}
\put(60,2){\line(1,0){50}}
\multiput(60,2)(10,0){6}{\circle*{1}}
\put(65,2){\oval(10,8)[t]}
\put(95,2){\oval(10,8)[t]}
\put(95,2){\oval(30,16)[t]}
\put(120,2){\line(1,0){50}}
\multiput(120,2)(10,0){6}{\circle*{1}}
\put(165,2){\oval(10,8)[t]}
\put(135,2){\oval(10,8)[t]}
\put(135,2){\oval(30,16)[t]}

\put(180,2){\line(1,0){50}}
\multiput(180,2)(10,0){6}{\circle*{1}}
\put(215,2){\oval(10,8)[t]}
\put(195,2){\oval(10,8)[t]}
\put(205,2){\oval(50,16)[t]}
\put(240,2){\line(1,0){50}}
\multiput(240,2)(10,0){6}{\circle*{1}}
\put(265,2){\oval(10,8)[t]}
\put(265,2){\oval(50,24)[t]}
\put(265,2){\oval(30,16)[t]}
\end{picture}\\ \\
\noi  The following result is known as the Fundamental Recurrence for Catalan numbers \cite{St16} 1.2.
\be \label{W1}C_{r+1}  = \sum_{i=1}^{r+1} C_{r-i+1}C_{i-1}. \ee  
In Theorem \ref{Lemma2}, 
we give a {\it Catalan-type recursion} for the set $\flat f$ from \eqref{VGL}.  
That is we show that $ \flat f$ is a disjoint union over an index set $I$ 
\be \label{W7}\flat f =\;\du_{i \in I}\flat f_i \ee where each $\flat f_i$ is (up to a natural bijection) a Cartesian product of sets
\be \label{W9} \flat f_i = \flat f_{L, i }\ti \flat  f_{R, i } .\ee
Then clearly 
\be \label{W}|\flat f|  = \sum_{i \in I}  |
\flat f_{L, i }
| |\flat  f_{R, i }|. \nn\ee    
\noi 
We postpone some of the definitions in Equations \eqref{W7} and \eqref{W9}, but the idea is that the set $\flat f$ can be enumerated by listing ``products'' of sets $\flat g \ti \flat h$ where $\# g, \# h < \# f$.
\br \label{R} {\rm 
Consider the problem of connecting $2n$ points lying on a  line by $n$ nonintersecting caps (as before), with the requirement that the leftmost point is connected to the rightmost by a cap $D$. The number of solutions to this problem is   $C_{n-1}$.  We call this the bounded definition  of Catalan numbers. (All other caps are bounded by $D$). }
\er 

\subsection{Weight and cap diagrams} \label{ab}
Since the Catalan number $C_n$ was defined  
without numbering the $2n$ points to be connected, it seems likely that the set $\flat f$ could be enumerated without the use of coordinates, using lattice paths. However we follow the tradition in representation theory of describing highest weights using their ``numerical marks."  The details are unimportant for our main result.  If suffices to say that there is a bijection between $X(r)$ and the set $F_r$
defined in the next paragraph, for details see \cite{MSe}.  
We use the same symbol for  $f\in X(r)$  and the corresponding element of $F_r$. The material on Kac modules was included solely to motivate the combinatorial problem and we have no more use for the set $X(r)$.
\\ \\
We use the language of weight and cap diagrams introduced by Brundan \cite{BS4}, \cite{MSe}. 
Let $F$  be
the set of all functions from $\mathbb{Z}$
 to the set $ \{ \ti, \ci\}$ such that $f(a) = \ci$ for all except finitely
 many  $a \in \mathbb{Z}$.  
If 
\be \label{lab}f^{-1}(\ti) =\{a_1 <a_2 \cdots <a_r\},\ee 
we also  write $f =(a_1 , a_2, \cdots , a_r)$ and set $\# f =r$,  $F_r=\{f\in F|\# f =r\}$. We call  $\# f $ the {\it degree of atypicality} of $f$.  We say $f$ and $g =(b_1 , b_2, \cdots , b_r)$ are equal  {\it up to shift} if $a_i - b_i$ is constant.  
\\ \\
The  {\it weight diagram}
$D_{wt}(f)$ is  a  number line with the symbol $\ti$ with 
 at each $a_i$ and symbols $\ci$ at all other integers.  
In Example \ref{xx} when $f=(1,2,3,7,9)$ we draw 
$D_{wt}(f)$ together with some other information.
\\ \\
A {\it cap} $C$ is the upper half of a circle joining two
integers $a$ and $b$. 
If $b < a$ we
say that $C$ {\it starts} at $b$ and {\it ends} at $a$ and we write
$b(C) = b$, and $e(C) = a.$ A finite set of caps  
 is called a {\it cap diagram}
if no two caps intersect, and no  integers remain inside the caps, unless they are ends of other caps.
\\ \\
\noindent If $f =(a_1 , a_2, \cdots , a_r)$, the  {\it cap diagram}    $D_{cap}(f)$ has a cap starting at each $a_i$. To find the ends of the
caps start with $a_r$ and work from right to left. The end of the
cap starting at $b$  is located  at the leftmost available symbol $\ci$ which is to the
right of $b$.
\\ \\
We say that a weight diagram and a cap diagram {\it match} if  when superimposed on the same number line,
each cap connects a $\ci$ to a $\ti$. 
In Section \ref{sec1}, given   $f \in F_r$ we explain how to enumerate the set
$\flat f$ from \eqref{VGL}. 
 We draw all cap diagrams match that match $f=(2,4)$ numbering the points where  the caps start.
\\ \\
\setlength\unitlength{0.48mm}
\begin{picture}(30,20)
\put(10,2){\line(1,0){40}}
\multiput(10,2)(10,0){5}{\circle*{1}}
\put(15,2){\oval(10,8)[t]}
\put(35,2){\oval(10,8)[t]}
\put(17,0){$\ti$}
\put(37,0){$\ti$}
\put(260,-12){\tiny{2}}
\put(200,-12){\tiny{2}}
\put(130,-12){\tiny{1}}
\put(70,-12){\tiny{1}}
\put(10,-12){\tiny{1}}
\put(270,-12){\tiny{3}}
\put(220,-12){\tiny{4}}
\put(140,-12){\tiny{2}}
\put(100,-12){\tiny{4}}
\put(30,-12){\tiny{3}}

\put(70,2){\line(1,0){40}}
\multiput(70,2)(10,0){5}{\circle*{1}}
\put(75,2){\oval(10,8)[t]}
\put(105,2){\oval(10,8)[t]}
\put(77,0){$\ti$}
\put(97,0){$\ti$}
\put(130,2){\line(1,0){40}}
\multiput(130,2)(10,0){5}{\circle*{1}}
\put(145,2){\oval(10,8)[t]}
\put(145,2){\oval(30,16)[t]}
\put(157,0){$\ti$}
\put(137,0){$\ti$}
\put(217,0){$\ti$}
\put(197,0){$\ti$}
\put(190,2){\line(1,0){40}}
\multiput(190,2)(10,0){5}{\circle*{1}}
\put(225,2){\oval(10,8)[t]}
\put(205,2){\oval(10,8)[t]}
\put(250,2){\line(1,0){40}}
\multiput(250,2)(10,0){5}{\circle*{1}}
\put(275,2){\oval(10,8)[t]}
\put(275,2){\oval(30,16)[t]}
\put(277,0){$\ti$}
\put(257,0){$\ti$}
\end{picture}
\\ \\ \\ Thus when $f= (2,4)$ we have  
\be \label{24}
\flat f = \{(1,3), \; (1, 4), \;  (1,2) , \;(2,4), \; (2,3)  \}.\ee
Our main result was conjectured in \cite{MSe}.
\bt \label{T 1} If $\# f =r$, then $|\flat f| \le C_{r+1}$ with equality iff $f =p$ up to shift, where  
$p = (2, 4, \ldots, 2r)$. \et
\noi A key step is to write $\flat f$ in the form \eqref{W7} and to prove a bound on the size of the index set $I = \LN f$ which will be defined in  Section \ref{sei}. 
\bp \label{Prop 1} If $\# f =r$, then $|\LN f| \le {r+1}$ with equality iff $f =p$ up to shift. \ep \noi 
This is proved using the following fact which is a Corollary to the Extreme Value Theorem. {\it If $\ttT:\R\lra \R$ is a continuous function that is not constant on any open interval, then any two zeroes of $\ttT$ are separated by a local extremum.}  
In fact our function $\ttT$ will have finitely many zeroes and  $\lim_{x\ra\pm \infty} \ttT(x) = \mp \infty$. This implies that $\ttT$ has an equal number of  local maxima and 
 minima.  {\it Also any two local maxima $($resp. 
 minima$)$ of $\ttT$ are separated by a local  
 minimum $($resp. local maximum$)$.}\\ \\
It is often necesssary  to count  the number of symbols $ \ti$ minus the number of symbols 
$\ci$  contained in the interval  $(b,a)$ of the weight diagram $D_{wt}(f)$. This statistic, which is related to tally functions in Lemma  \ref {Lemma3},
 is used in some form 
in various papers: \cite{GS},  \cite{GS2}, \cite{GH},  \cite{MSe} and \cite{SZ}.

\section{Tally functions and a Catalan-type recursion for   $\flat f$}

\subsection{Tallies} Corresponding to  each $f$ as in \eqref{lab} we define a  {\it tally function} $\ttT_f: \Z \lra \Z $.  The idea is to associate a non-negative integer (a tally) to each point on the number line so that, moving left to right, the tally increases (resp. decreases) by one each time a symbol $\ti$ (resp. $\ci$) is passed.  
Thus we require 
\be \label{nll} \ttT_f (b+1) =
\left\{ \begin{array}
  {cc} \ttT_f (b) +1&\mbox{if} \;\; f( b+1)= \ti
\\ \ttT_f (b) -1 &\mbox{ otherwise}
\end{array} \right.\nn\ee This only specifies $\ttT_f$ up to an additive constant, so we also assume 
\be \label{ml}\ttT_f(a) =1,\ee where $a=a_r$. This assumption gives  good window in which to view the graph of $\ttT_f$. It also relates potential moves to the zeros of  $ \ttT_f$. There is a unique way to extend 
$ \ttT_f$ 
to a function $\ttT_f:\R\lra \R$  by connecting the values at consecutive integers by line segments with slopes $\pm 1$. 
If  $b<a$ are integers set $[b,a ] = \{z\in \Z| b\le z\le a\}$.  We denote the corresponding real closed interval by $[b,a ]_\R$. 
Likewise for other kinds of intervals with at least one integer endpoint. 
\\ \\
For   $f$ as in \eqref{lab}, set $\I =[a_1 -1, a_r]. $ 
There is a disjoint union 
$$\R=(-\infty, a_1-1)_\R\cup \I_\R\cup (a_r, \infty)_\R.$$
The graph of $ \ttT_f$ consisits of two half-lines with slope $-1$ covering the first  and last intervals above, joined by a 
lattice path  $\lp(f)$  on $\I_\R$, 
with steps $(1,1)$ and $(1,-1)$.  
\bexa \label{xx} {\rm Here is the graph of $\ttT_f$ for $f=(1,2,3,7,9)$.  The lattice path 
 occurs between the arrows, and we give the values of the tallies here at integer points. The entries in $f$ are the location of the symbols $\ti$. 
\\ \\ \\ \\
\setlength\unitlength{0.48mm}
\begin{picture}(30,20)(-88,0)
\put(-20,2){\line(1,0){140}}
\put(10,2){\line(1,1){30}}
\put(70,2){\line(-1,1){30}}
\put(90,2){\line(-1,1){10}}
\put(90,2){\line(1,1){10}}
\put(70,2){\line(1,1){10}}
\put(10,2){\vector(-1,1){28}}
\put(100,12){\vector(1,-1){20}}
\multiput(-20,2)(10,0){14}{\circle*{1}}
\multiput(10,2)(10,0){10}{\circle*{1}}
\put(37,0){$\ti$}
\put(8.5,-10){\tiny{0}}
\put(18.5,-10){\tiny{1}}
\put(28.5,-10){\tiny{2}}
\put(38.5,-10){\tiny{3}}
\put(48.5,-10){\tiny{2}}
\put(58.5,-10){\tiny{1}}
\put(68.5,-10){\tiny{0}}
\put(78.5,-10){\tiny{1}}
\put(88.5,-10){\tiny{0}}
\put(98.5,-10){\tiny{1}}
\put(17,0){$\ti$}
\put(37,0){$\ti$}
\put(27,0){$\ti$}
\multiput(70,2)(10,0){05}{\circle*{1}}
\put(77,0){$\ti$}
\put(97,0){$\ti$}
\end{picture}}
\eexa
\vspace{0.4cm}

\noi 
A diagram of the form \\\ 
\setlength\unitlength{0.48mm}
\begin{picture}(30,20)(-88,0)
\put(60,2){$\cdots$}
\put(10,2){\line(1,1){10}}
\put(30,2){\line(-1,1){10}}
\put(30,2){\line(1,1){10}}
\put(50,2){\line(-1,1){10}}
\put(90,2){\line(-1,1){10}}
\put(90,2){\line(1,1){10}}
\put(10,2){\vector(-1,1){10}}
\put(100,12){\vector(1,-1){10}}
\multiput(-10,2)(10,0){6}{\circle*{1}}
\put(-10,2){\line(1,0){65}}
\multiput(80,2)(10,0){6}{\circle*{1}}
\put(80,2){\line(1,0){45}}
\end{picture}\\ with $r$ maximum values of 1, and $r$ minimum values of $0$ is called a {\it zigzag diagram of rank $r$.}

 \subsection{Moves: potential and legal} \label{sei} 
\bl \label{Lemma3}
If $b<a$ are integers, then the number of symbols $ \ti$ minus the number of symbols 
$\ci$  contained in the interval  $(b,a)$ of the weight diagram $D_{wt}(f)$ equals  $\ttT_f (a-1)- \ttT_f (b).$ 
\el 
\bpf Simple counting or induction on $a$.\epf
\noi
If $b<a$, $f(b)=\ci$ and $f(a)=\ti$, 
As in \cite{MSe}, define $f_{b} \in F,$ (resp. $f^a \in
F$) by adding $b$ to $\ti(f)$ (resp. deleting $a$ from $\ti(f)$).
We also set $f_{b}^a =
(f_{b})^a =  (f^a)_{b}.$\\
\\ 
We say that the interval $ (b,a)$ is {\it balanced} if it contains an equal number of symbols $\ti$ and $\ci$. If  $ (b,a)$ is { balanced} and $f(b)=\ci$
we say that $f\lra f^a_b $ is a {\it potential move}. 
If there is a potential move $f\lra f^a_b $, then for some integer $i\in [r]$, there is an equal number $i-1$ of symbols $\ti$ and $\ci$ in the interval $(b,a)$ and we have 
  \be \label{zt} a-b =2i-1.\ee 
\noi Set 
\be \label{zf} \PM f=\{i\in [r]|  \mbox{ there is a  potential move  } f \lra f^a_b   \mbox{ where    }  b=a+1-2i\}.\ee
Now suppose $i\in \PM f$ and set $b=a+1-2i$.  If in addition  

\be \label{z9} \ttT_f (x-1)- \ttT_f (b) \ge 0 \mbox{ for all } x\in (b,a]\ee
 and $\ttT_f (a-1)- \ttT_f (b) =0$, we say $\L^a_b $ is  a {\it legal move } from $f$ to $g$ or there  is  a {\it legal move }
$\L^a_b f$. In  \cite{MSe} this is called a legal move of weight zero, but we have no need for legal  moves of positive weight. Set

\be \label{zf} \LM f=\{i\in [r]|  \mbox{ there is a  legal move  } \L^a_b f  \mbox{ where    }  b=a+1-2i\}.\nn\ee
By assumption \eqref{ml} the 
set of zeroes $\ttT_f^{-1}(0)$  of $\ttT_f$ satisfies  \be \label{z2} \{a+1-2i| i\in \PM f\} =
\{b \in \ttT_f^{-1}(0)| b\le a, f(b) =  \ci,  (b,a) \mbox{ balanced} \}.\ee
If $a=b$ in \eqref{zt}, then $i = 1/2$. 
Set $\LN f = \PM f\cup 
\{1/2\}$ and $\L^a_a f =f$. \\ \\
In Example \ref{xx},  $
\PM f =\{1, 2, 5\}$  with corresponding legal moves  $\L^{9}_8, \L^{9}_6, $ and $\L^{9}_0$.
Here is an example where not all potential moves are legal.
\bexa \label{xx1} {\rm 
For example if  
$f=(1,4,5,7)$, then $\PM =\{0, 2,  6\}$ and  $\PM =\{6\}$, that is the only legal move starting at 7 is $\L^{7}_6$.  Note that 
 $f\lra f^{7}_0$ and  
$f\lra f^{7}_2$ are not   legal moves because $\ttT_f(3)<0.$
}\eexa

\subsection{A recursion for   $\flat f$} \label{sec1} 
Suppose $f$ satisfies \eqref{lab}, set $a=a_r$ and recall the definition of $ \flat f$ from  \eqref{W7}. 
  Let $\bar f \in F_{r-1}$ be given by $\bar f^{-1}(\ti) =(a_1, \ldots, a_{r-1})$, and set
\[ \flat\bar f= \{\;g \in F_{r-1} \;| \; D_{cap}(g) \; \mbox{matches} \; D_{wt}(\bar f)  \}, \]
\[ \flat_{1/2} f = \{\;g \in \flat f \;| \; D_{cap}(g) \; \mbox{has a cap joining} \;  a
\; \mbox{to} \; a +1 \}. \]
Our recursion expands on ideas from \cite{MSe}  Section 5.  In particular the following result is \cite{MSe} Lemma 5.1.
 \bl \label{Lemma 1} There is a bijection
$\flat\bar f \longrightarrow \flat_{1/2} f$ such that $g =(b_1, \ldots , b_{r-1}) \in \flat\bar f$ 
maps to 
$h = (b_1, \ldots, b_{r-1}, a)$. Hence by induction $|\flat_{1/2} f| \le C_r$.\el \noi 
Given $f$
as above, we can assume by induction that we have found
$\flat\bar f$ and hence $\flat_{1/2} f\subset \flat f$. The set $\flat_{1/2}$ is the set of matching cap diagrams with  a cap $C$ joining $a$ to $a +1$.  The remainder are found by replacing $C$ by a cap $D$ ending  at $a$.  That is 
for each  $h \in \flat_{1/2} f$ we need list all weights obtained from $h$ by a legal move.   
 To do this,  for all $i\in \PM f$, set

\be \label{z71}  \flat_{i} f = \{\;g \in \flat f \;| \; D_{cap}(g) \; \mbox{has a cap joining} \;  a
\; \mbox{to} \; b=a +1-2i \}. \ee
 Define new weights $f_{L, i }, f_{R, i } \in F$ by  
\be \label{z21}  f_{L, i }^{-1}(\ti) =f^{-1}(\ti) \cap [a_1, b), \quad 
 f_{R, i }^{-1}(\ti) =f^{-1}(\ti)  \cap (b, a_r].
\ee 
\bl There is a bijection
\be \label{V9} \flat f_i  \llra \flat f_{L, i }\ti \flat  f_{R, i }.\ee
\el
\bpf This  follows since $\flat f_{L, i }$ $($resp. $\flat  f_{R, i })$  is the set of cap diagrams that can be placed with start to the left of $D$  $($resp. under $D)$   resulting in a diagram in $\flat f_i$.
\epf
\bl
 If $g\in  \flat_{i} f$ and suppose  $b = a+1-2i$.  Then $h= g^b_a\in  \flat_{1/2} f$ and there is a legal move $h\lra g$. 
\el
\bpf 
 Let $D$ be a cap joining $b$ to $a$.  
Let $S$ be the set of cap diagrams that can be drawn under $D$ without regarding the location of the  $n=i-1$ symbol pairs in 
$D_{cap}(g) $. 
By Remark \ref{R}, $|S| = C_{i-1}$.
For each such diagram replace the left (resp. right) end of each cap by the number 1 (resp. $-1$). We obtain 
a ballot sequence (as in \cite{St16} Section 1.5) of length  $2n$  consisting of $1$'s and $-1$'s such that each partial sum is nonnegative.  This applies to the portion of  
$D_{cap}(g) $ covering the interval $[b, a]$ and is precisely the condition \eqref{z9} required for a legal move $h\lra g$.
\epf
\bc We have $$\flat f= \{ g \in F|  \mbox{there is a legal move } h\lra g \mbox{ for some } h\in  \flat_{1/2} f\}.$$ 
\ec

\bt \label{Lemma2} There is a {\it Catalan-type recursion} for the set $\flat f$ as in Equations  \eqref{W7} and \eqref{W9}, where
$\flat_{i} f$,  and 
$f_{L, i }, f_{R, i } \in F$ are given by  
Equations \eqref{z71} and \eqref{z21} respectively.
\et \noi
\bpf 
By construction each set $\flat_{i} f$ is contained in $\flat f$ and $\flat_{i} f \cap \flat_{j} f=\emptyset$  if $i\neq j$,.  Also if $g\in \flat f  \backslash  \flat_{1/2} f$, then $g$ has a cap joining $a$ to $ a +1-2i$ for some $i$, so $g \in \flat_{i} f$. Thus \eqref{W7} holds and  \eqref{W9} holds by \eqref{V9}.
\epf
\noi There is a restriction on  $ \flat f_{R, i }$ because all caps that start  under $D$ must also end under $D$.    There is no corresponding restriction on $\flat f_{L, i }$: a cap that starts to the left of $D$ may end to the left or right of $D$. 
We have $\# f_{L, i }={r-i}=j$ and  $\#   f_{R, i }  =i-1$ Hence by induction and 
Remark \ref{R} we obtain
  
\bl \label{Lemmaz} There are  upper bounds
$|\flat f_{L, i }| \le C_{j+1} =C_{r+1-i}$ and $|\flat  f_{R, i }|\le C_{i-1}.$  
\el 
\noi 
\bexa{\rm 
 If 
$p = (2, 4, \ldots, 2r)$, then $|\flat p| =C_{r+1}$.  This is  shown in \cite{MSe} Example 2.4.   To prove the Fundamental Recurrence \eqref{W1},  note that 
 $\PM p =[r]$,  and for $i\in [r]$, $|\flat  p_{L, i }|  = C_{r+1-i}$ , $|\flat  p_{R, i }| = C_{i-1}$.  Then use  Theorem \ref{Lemma2}. 
}\eexa

\bexa \label{xx1} {\rm 
This example shows that it is necessary to use potential, rather than legal moves  in Lemma \ref{Lemma2}.  Let 
$f=(2,3), g= (1,3), h=(0,1)$.   There are no legal moves starting with 3, but   $\PM f =\{2\}$.  We have $ \flat_{1/2} f =\{f,g\}$, $\flat_{2}=\{h\}$
and $ \flat f$ is the disjoint union of these two sets}\eexa

\subsection{Local behaviour of tally functions} \label{sec3} 
When the two sided derivative of $\ttT_f  $ exists at $b$, we denote it by $\ttT'_f(b) $. 
Now $\Z$ is the disjoint union of the four sets  
$$
\max f =  \{b \in \Z| \ttT_f(b) \mbox{ has a local maximum at  } \;b\},$$ $$ 
  \min f = \{b  \in \Z| \ttT_f(b) \mbox{ has a local minimum at  } \;b\},$$
 $$f_{\pm} = \{b\in \Z| \ttT'_f(b)=\pm1\}.$$

\bl \label{LE1}
There are disjoint unions 
\be \label{EQ1} 
f^{-1}(\ti)= {\max} f \cup f_{+}   
\ee
\be \label{EQ2} 
f^{-1}(\ci)\cap \I = {\min} f \cup (f_{-} \cap \I)  
\ee
\el \bpf This follows by looking at 
the diagrams below, each of which corresponds to one of the sets in the disjoint union. The symbol $*$  denotes either  $\ci$ or  $\ti$. 
\noi \\ \\ \\ \\
\setlength\unitlength{0.58mm}
\begin{picture}(30,20)
\put(10,2){\line(1,0){40}}
\put(43,2){\circle*{1}}
\put(27,0){$\ti$}\put(27,-10){$b$}
\put(14,0){$\large{ *}$}
\put(70,2){\line(1,0){40}}
\put(74,0){$\large{ *}$}
\put(87.5,0){$\ti$} \put(87.5,-10){$b$}
\put(101,0){$\ti$}
\put(147.5,15){\line(-1,1){13}}
\put(147.5,15){\line(1,1){13}}
\put(77,10){\line(1,1){26}}
\put(17,15){\line(1,1){13}}
\put(43,15){\line(-1,1){13}}
\put(223,10){\line(-1,1){26}}
\put(209.5,2){\circle*{1}}
\put(149.5,2){\circle*{1}}\put(148,-10){$b$} \put(208.5,-10){$b$}
\put(130,2){\line(1,0){40}}
\put(161,0){$\ti$}
\put(134,0){$\large{ *}$}
\put(223,2){\circle*{1}}
\put(194,0){$\large{ *}$}
\put(190,2){\line(1,0){40}}
\end{picture}\epf

\subsection{Proof of the main results} \label{
bookkeeping} First we handle some  special cases.  A subset $A$ of $\Z$ is {\it connected} if $i,j\in A, i<j$ implies $[i,j] \subset A.$
\bl \label{LL1} If $\# f =r$, then 
$|\PM f|\le r-1$ provided one of the following holds.
\bi \itema $f_{+} \neq\emptyset$ 
\itemb $f_+ =\emptyset$ and $f_-\neq\emptyset$ 
\ei
\el
\bpf (a) Set $S= \{b \in \ttT_f^{-1}(0)| b\le a, f(b) =  \ci\}$.  By assumption and \eqref{EQ1}, $|\max f| = |\min f| \le r-1$.   
Then $\PM f= |S|$.  Write $f^{-1}(\ci)\cap \I $ as a disjoint union 
$$f^{-1}(\ci)\cap \I =\du_{i=1}^k A_i$$ of connected subsets $A_i$ with $k$ minimal.  Then $k=r-1$  because there is one local max between two consecutive such sets and another  to the right of all $A_i$.  Since  $|S\cap A_i| \le 1$ for each $i$,  we obtain  $|S|\le r-1.$ \\ \\ 
(b) Suppose $f_+ =\emptyset$  and $c\in \I$ is maximal such that $\ttT'_f(c)=-1$ and $a_{j}<c<a_{j+1}$.
Then the part of the  graph  of $\ttT_f $  
covering $[a_{j+1}  -1, a_r]$ is a zigzag diagram of rank $r-j$.  Also for $b<c$ we have $\ttT_f (b) \ge \ttT_f(a). $ 
 Thus the only legal moves  
$\L^a_b f$ have $b =a+1-2i$, where $i\in [r-j]$.\epf
\noi {\it Proof of {Proposition} \ref{Prop 1}}  
Now assume 
$f_+=f_- =\emptyset$. 
 The graph  of $\ttT_f $  is a zigzag diagram of rank $r$. 
By shifting, we can assume that $a=2r$ is the rightmost entry in 
$f^{-1}(\ti)$   and then it follows that $f=p,$
\hfill $\Box$
\\ \\ 
\noi {\it Proof of {Theorem} \ref{T 1}} We use Lemmas \ref{Lemma 1}, Theorem \ref{Lemma2},  Lemma \ref{Lemmaz} and then Equation \eqref{W1} 
\by \label{b11}
|\flat f|& =&\sum_{i \in \LN f}|\flat_i  f| \nn\\
&=&|\flat_{1/2} f| + \sum_{i \in \PM f}|
 \flat  f_{L, i }| | \flat f_{R, i }|\nn \\
 &\le& C_0 C_r+ \sum_{i \in \PM f\subseteq [r]}C_{r-i+1}C_{i-1}\nn\\
&\le& \sum_{i=1}^{r+1} C_{r-i+1}C_{i-1} = C_{r+1}.\nn\ey 
If equality holds, then $\PM f= [r]$, so $f =p$ up to  shift
by Proposition \ref{Prop 1}.\hfill $\Box$
\br{\rm 
If $\# f =r$ and   $\ttT_f$ has more than one local minimum, then $|\flat f| > {r+1}$.  This can be used to show that  $|\flat f| \ge {r+1}$ with equality iff $f =q$ up to shift, where  
$q = (1, 2, \ldots, r)$.  The Kac module $K(q)$ is a {\it Kostant module} as in the introduction to 
\cite{BS4} or \cite{GH} Definition 3.5.3.  Such a module has a multiplicity-free resolution by Kac modules.
}\er
\bexa {\rm We give an example to show how to enumerate $\flat f$ in practice.  Let $f=(2, 3, 5, 7)$, so $a=7$.  This weight can be obtained from \cite{SZ} (5.3)   by removing all core symbols.  In the table below we list all $g$ such that $D_{cap}(g) $ matches $D_{wt}(f)$. The top row  lists all $g\in  \flat_{1/2} f$.  Consider the column headed by $g$.  If $g(b)=\ti$ then $g \lra M^7_b g$ 
is not even a potential move. The entry corresponding to  $M^7_b g$ is blank in this case.   
For each  $b\in \PM = \{0,4,6\}$, the move   $M^7_b$ is listed and the result in this row is  $M^7_a g$ provided this is in $\flat f.$
\[ \begin{tabular}{|c||c|c|c|c|c|c|c|c|} \hline
& 2357 & 2347&1357&0157&0137&1347&0147\\ \hline \hline 
 $M^7_6$
  & 2356 & 2346&1356&0156&0136&1346&0146\\
  \hline
$M^7_4$ & 2345 &  &1345&0145&0134& &\\
    \hline
$M^7_0$ &{\color{red}  0235} & {\color{red}  0234}&0135&&&{\color{red}  0134}&\\
    \hline
\end{tabular}\]
If $h$ is one of the entries in red, 
then $g\lra h$ is not a legal move. Note that in this abbreviated notaion $0235$ and  0234 do not belong to $\flat f$.  However $0134 \in \flat f$ because $0137 \lra 0134$ is a legal move. 
Ignoring the entries in red we obtain the $19 $ highest weights of the composition factors of $K(f)$.

}\eexa


\end{document}